\documentclass[11pt,reqno,oneside]{amsproc}
\title[Remarks on the BKM criterion]{Remarks on the BKM criterion for both classical and variable Euler systems}

\author[B.~Ingimarson]{Benjamin Ingimarson}
\address{Department of Mathematics, University of Southern California, Los Angeles, CA 90089}
\email{ingimars@usc.edu}

  \chardef\forshowkeys=0
  \chardef\showllabel=0
  \chardef\refcheck=0
  \chardef\sketches=0
  \chardef\figures=1

\usepackage{enumitem,mathtools}
\usepackage{datetime}
\usepackage{fancyhdr}
\usepackage{comment}
\usepackage{listings}

\allowdisplaybreaks
\ifnum\forshowkeys=1

  \usepackage[notref,notcite,color]{showkeys}
\fi
\usepackage[margin=1in]{geometry}
\usepackage{amsmath, amsthm, amssymb}
\usepackage{times}
\usepackage{graphicx}
\usepackage[usenames,dvipsnames,svgnames,table]{xcolor}
\usepackage{marginnote}
\usepackage[unicode,breaklinks=true,colorlinks=true,linkcolor=blue,urlcolor=blue,citecolor=blue]{hyperref}
\usepackage[most]{tcolorbox}
\usepackage{tikz}

\ifnum\refcheck=1
  \usepackage{refcheck}
\fi

\usepackage{import}
\usepackage{xifthen}
\usepackage{pdfpages}
\usepackage{transparent}

\begin{document}
\def\YY{X}
\def\OO{\mathcal O}
\def\SS{\mathbb S}
\def\CC{\mathbb C}
\def\RR{\mathbb R}
\def\TT{\mathbb T}
\def\ZZ{\mathbb Z}
\def\HH{\mathbb H}
\def\RSZ{\mathcal R}
\def\LL{\mathcal L}
\def\SL{\LL^1}
\def\ZL{\LL^\infty}
\def\GG{\mathcal G}
\def\tt{\langle t\rangle}
\def\erf{\mathrm{Erf}}
\def\mgt#1{\textcolor{magenta}{#1}}
\def\ff{\rho}
\def\gg{G}
\def\sqrtnu{\sqrt{\nu}}
\def\ww{w}
\def\ft#1{#1_\xi}
\def\les{\lesssim}
\def\lec{\lesssim}
\def\ges{\gtrsim}
\def\gec{\gtrsim}
\renewcommand*{\Re}{\ensuremath{\mathrm{{\mathbb R}e\,}}}
\renewcommand*{\Im}{\ensuremath{\mathrm{{\mathbb I}m\,}}}
\ifnum\showllabel=1
 \def\llabel#1{\marginnote{\color{lightgray}\rm\small(#1)}[-0.0cm]\notag}
%%%%%%%%%%%%%%%%%%%\def\llabel#1{\label{#1}}
%  \reversemarginpar
  %\def\llabel#1{\notag}
%  \def\llabel#1{\label{#1}}
\else
% % \def\llabel#1{\nonumber}
 \def\llabel#1{\notag}
\fi
\newcommand{\norm}[1]{\left\|#1\right\|}
\newcommand{\nnorm}[1]{\lVert #1\rVert}
\newcommand{\abs}[1]{\left|#1\right|}
\newcommand{\NORM}[1]{|\!|\!| #1|\!|\!|}

\newtheorem{Theorem}{Theorem}[section]
\newtheorem{Corollary}[Theorem]{Corollary}
\newtheorem{Proposition}[Theorem]{Proposition}
\newtheorem{Lemma}[Theorem]{Lemma}
\newtheorem{Remark}[Theorem]{Remark}
\newtheorem{definition}{Definition}[section]

\def\theequation{\thesection.\arabic{equation}}
\numberwithin{equation}{section}
%%%%%%%%%%%%%%%%%%%%%%%%%%
\definecolor{mygray}{rgb}{.6,.6,.6}
\definecolor{myblue}{rgb}{9, 0, 1}
\definecolor{colorforkeys}{rgb}{1.0,0.0,0.0}
\newlength\mytemplen
\newsavebox\mytempbox
\makeatletter
\newcommand\mybluebox{%
    \@ifnextchar[%]
       {\@mybluebox}%
       {\@mybluebox[0pt]}}
\def\@mybluebox[#1]{%
    \@ifnextchar[%]
       {\@@mybluebox[#1]}%
       {\@@mybluebox[#1][0pt]}}
\def\@@mybluebox[#1][#2]#3{
    \sbox\mytempbox{#3}%
    \mytemplen\ht\mytempbox
    \advance\mytemplen #1\relax
    \ht\mytempbox\mytemplen
    \mytemplen\dp\mytempbox
    \advance\mytemplen #2\relax
    \dp\mytempbox\mytemplen
    \colorbox{myblue}{\hspace{1em}\usebox{\mytempbox}\hspace{1em}}}
\makeatother
%%%%%%%%%%%%%%%%%%%%%%%%%%

%Igor's macros  varmac
\def\BMO{\text{BMO}}
\def\pv{\text{p.v.}}
\def\bnew{\color{red}}
\def\enew{\color{black}}
\def\bold{\color{blue}}
\def\eold{\color{black}}
\def\rr{r}
\def\weaks{\text{\,\,\,\,\,\,weakly-* in }}
\def\inn{\text{\,\,\,\,\,\,in }}
\def\cof{\mathop{\rm cof\,}\nolimits}
\def\Dn{\frac{\partial}{\partial N}}
\def\Dnn#1{\frac{\partial #1}{\partial N}}
\def\tdb{\tilde{b}}
\def\tda{b}
\def\qqq{u}
\def\lat{\Delta_2}
\def\biglinem{\vskip0.5truecm\par==========================\par\vskip0.5truecm}
\def\inon#1{\hbox{\ \ \ \ \ \ \ }\hbox{#1}}                %in or on
\def\onon#1{\inon{on~$#1$}}
\def\inin#1{\inon{in~$#1$}}
\def\FF{F}
\def\andand{\text{\indeq and\indeq}}
\def\ww{w(y)}
\def\ll{{\color{red}\ell}}
\def\ee{\mathrm{e}}
\def\startnewsection#1#2{ \section{#1}\label{#2}\setcounter{equation}{0}}   
%\starts a new section\def\dist{\mathop{\rm dist}\nolimits}
\def\nnewpage{ }
\def\sgn{\mathop{\rm sgn\,}\nolimits}    
\def\Tr{\mathop{\rm Tr}\nolimits}    
\def\div{\mathop{\rm div}\nolimits}
\def\curl{\mathop{\rm curl}\nolimits}
\def\dist{\mathop{\rm dist}\nolimits}
\def\id{\mathop{\rm id}\nolimits}
\def\supp{\mathop{\rm supp}\nolimits}
\def\indeq{\quad{}}           
\def\period{.}                       
\def\semicolon{\,;}                  
\def\cmi#1{\text{~{{\coli IK: \underline{#1}}}~}}
\def\coli{\color{colorigor}}
\definecolor{colorigor}{rgb}{.5, 0.2, 0.8}
\def\colr{\color{red}}
\def\colrr{\color{black}}
\def\colb{\color{black}}
\def\coly{\color{lightgray}}
\definecolor{colorgggg}{rgb}{0.1,0.5,0.3}
\definecolor{colorllll}{rgb}{0.0,0.7,0.0}
\definecolor{colorhhhh}{rgb}{0.3,0.75,0.4}
\definecolor{colorpppp}{rgb}{0.7,0.0,0.2}
%\definecolor{coloroooo}{rgb}{0.9,0.4,0}
\definecolor{coloroooo}{rgb}{0.45,0.0,0.0}
\definecolor{colorqqqq}{rgb}{0.1,0.7,0}
\def\colg{\color{colorgggg}}
\def\collg{\color{colorllll}}
\def\cole{\color{coloroooo}}
\def\coleo{\color{colorpppp}}
\def\colu{\color{blue}}
\def\colc{\color{colorhhhh}}
\def\colW{\colb}   %color for weight
\definecolor{coloraaaa}{rgb}{0.6,0.6,0.6}%%%out
\def\colw{\color{coloraaaa}}
\def\comma{ {\rm ,\qquad{}} }            
\def\commaone{ {\rm ,\quad{}} }          
\def\les{\lesssim}
\def\nts#1{{\color{red}\hbox{\bf ~#1~}}} 
\def\ntsf#1{\footnote{\color{colorgggg}\hbox{#1}}} 
\def\blackdot{{\color{red}{\hskip-.0truecm\rule[-1mm]{4mm}{4mm}\hskip.2truecm}}\hskip-.3truecm}
\def\bluedot{{\color{blue}{\hskip-.0truecm\rule[-1mm]{4mm}{4mm}\hskip.2truecm}}\hskip-.3truecm}
\def\purpledot{{\color{colorpppp}{\hskip-.0truecm\rule[-1mm]{4mm}{4mm}\hskip.2truecm}}\hskip-.3truecm}
\def\greendot{{\color{colorgggg}{\hskip-.0truecm\rule[-1mm]{4mm}{4mm}\hskip.2truecm}}\hskip-.3truecm}
\def\cyandot{{\color{cyan}{\hskip-.0truecm\rule[-1mm]{4mm}{4mm}\hskip.2truecm}}\hskip-.3truecm}
\def\reddot{{\color{red}{\hskip-.0truecm\rule[-1mm]{4mm}{4mm}\hskip.2truecm}}\hskip-.3truecm}
\def\tdot{{\color{green}{\hskip-.0truecm\rule[-.5mm]{3mm}{3mm}\hskip.2truecm}}\hskip-.1truecm}
\def\gdot{\greendot}
\def\bdot{\bluedot}
\def\ydot{\cyandot}
\def\rdot{\cyandot}
\def\fractext#1#2{{#1}/{#2}}
\def\ii{\hat\imath}
%%%%%%%
\def\fei#1{\textcolor{blue}{#1}}
\def\vlad#1{\textcolor{cyan}{#1}}
\def\igor#1{\text{{\textcolor{colorqqqq}{#1}}}}
\def\igorf#1{\footnote{\text{{\textcolor{colorqqqq}{#1}}}}}
%varmac
%\def\AA{Y}
\def\Omf{\Omega_{\text f}}
\def\Ome{\Omega_{\text e}}
\def\Omb{\Omega_{\text b}}
\def\Gaf{\Gamma_{\text f}}
\def\Gae{\Gamma_{\text e}}
\def\Gab{\Gamma_{\text b}}
\def\Gac{\Gamma_{\text c}}
\def\Nf{N^{\text f}}
\def\Ne{N^{\text e}}

\newcommand{\p}{\partial}
\renewcommand{\d}{\mathrm{d}}
\newcommand{\UE}{U^{\rm E}}
\newcommand{\PE}{P^{\rm E}}
\newcommand{\KP}{K_{\rm P}}
\newcommand{\uNS}{u^{\rm NS}}
\newcommand{\vNS}{v^{\rm NS}}
\newcommand{\pNS}{p^{\rm NS}}
\newcommand{\omegaNS}{\omega^{\rm NS}}
\newcommand{\uE}{u^{\rm E}}
\newcommand{\vE}{v^{\rm E}}
\newcommand{\pE}{p^{\rm E}}
\newcommand{\omegaE}{\omega^{\rm E}}
\newcommand{\ua}{u_{\rm   a}}
\newcommand{\va}{v_{\rm   a}}
\newcommand{\omegaa}{\omega_{\rm   a}}
\newcommand{\ue}{u_{\rm   e}}
\newcommand{\ve}{v_{\rm   e}}
\newcommand{\omegae}{\omega_{\rm e}}
\newcommand{\omegaeic}{\omega_{{\rm e}0}}
\newcommand{\ueic}{u_{{\rm   e}0}}
\newcommand{\veic}{v_{{\rm   e}0}}
\newcommand{\up}{u^{\rm P}}
\newcommand{\vp}{v^{\rm P}}
\newcommand{\tup}{{\tilde u}^{\rm P}}
\newcommand{\bvp}{{\bar v}^{\rm P}}
\newcommand{\omegap}{\omega^{\rm P}}
\newcommand{\tomegap}{\tilde \omega^{\rm P}}
\newcommand{\eps}{\varepsilon}  %ben's macro
\newcommand{\eqnb}{\begin{equation}}
\newcommand{\eqne}{\end{equation}}
  
\renewcommand{\up}{u^{\rm P}}
\renewcommand{\vp}{v^{\rm P}}
\renewcommand{\omegap}{\Omega^{\rm P}}
\renewcommand{\tomegap}{\omega^{\rm P}}

\begin{abstract}
In this paper, we establish new results on BKM-type criteria for the incompressible Euler equations on $\mathbb{R}^3$ as well as the Euler equations with space-time variable coefficients in a bounded domain. 
Our first result shows that if a finite-energy solution to the Euler equations on $\mathbb{R}^3$ blows up at a time $T_\ast$, and $T_\ast$ is the first such time, then $\limsup_{t\to T_\ast^-}(T_\ast -t)^{1+2k/5}\|D^k\omega(t)\|_{L^\infty} = \infty$ for every $k \geq 1$, thereby excluding the endpoint rate allowed by previously established lower bounds. 
The proof assumes the endpoint rate as an upper bound near the blowup time and uses a result by Chae and Wolf to obtain a contradiction. 
The second result gives a BKM-type criterion for solutions to the incompressible Euler equations with variable coefficients in the minimal regularity class $H^r$, where $r \in (2.5,3)$.
\end{abstract}

\maketitle

\setcounter{tocdepth}{2} 
\section{Introduction}\label{sec00}
First, we consider the 3D incompressible Euler equations,
  \begin{align}
   \begin{split}
    u_t
    + u\cdot \nabla u
    + \nabla p
    &= 0
    ,
    \\
    \div u &= 0,
   \end{split}
   \label{EQ00}
  \end{align}
in $\mathbb{R}^3 \times [0,T]$.
For an initial velocity field $u_0 \in H^r$, with $r > 2.5$, where $\|u_0\|_{H^r} \leq M$ for some $M > 0$, there exists a time $T = T(M) > 0$ such that~\eqref{EQ00} has a solution satisfying
  \begin{equation}
          u \in C([0,T];H^r(\mathbb{R}^3)) \cap C^1([0,T];H^{r-1}(\mathbb{R}^3))
   .
   \label{EQ01}
  \end{equation}
  It was shown by Beale, Kato, and Majda~\cite{BKM} that if such a solution exits the class~\eqref{EQ01} at a time $T_\ast$ and that $T_\ast$ is the first such time, then 
  \begin{equation}
   \int_0^{T_\ast} \|\omega(t)\|_{L^\infty}\,dt = \infty
   .
   \label{EQ01a}
  \end{equation}
Later works, such as~\cite{F,SY,Z}, proved~\eqref{EQ01a} in the case of a bounded domain. 
For improvements on~\eqref{EQ01a} by replacement with a weaker norm, see~\cite{A,C,CFM,P,KOT,KT},
and for a general exposition on the BKM criterion, see~\cite{MB,MP}.
In the case that $\Omega = \mathbb{R}^3$, there have been several recent advancements in finite-time blow-up of the Euler equations, see~\cite{Ch,CMZ,E,JMO,S,SWZZ}.
Regarding the structure of a putative singularity, Chae and Wolf~\cite{CW} showed that if a finite-energy solution on $\mathbb{R}^3$ satisfies the Type~I condition $\sup_{t<T_\ast}(T_\ast -t)\|\nabla u(t)\|_{L^\infty} < \infty$, then the kinetic energy cannot concentrate at isolated points at the time $T_\ast$; 
as a corollary, they excluded discretely self-similar blow-up in the energy-conserving scale under the same condition.

Under the assumption of blow-up, it is possible to derive lower bounds on the rate of growth of various norms, see~\cite{C,CIV,CP}. 
However, lower bounds on the maximum norm of vorticity were only recently shown in~\cite{IK1,IK2}.
In~\cite{IK1}, we established a lower bound on the blow-up rate of $\int_0^t \|\omega(s)\|_{L^\infty(\Omega)}\, ds$ when $\Omega = \mathbb{R}^3$, $\mathbb{T}^3$, or a smooth bounded domain, as well as pointwise-in-time lower bounds on $\|D^k\omega\|_{L^\infty(\Omega)}$ when $\Omega = \mathbb{R}^3$ or $\mathbb{T}^3$.
In~\cite{IK2}, we expanded the results on $\|D^k\omega\|_{L^\infty(\Omega)}$ to a smooth bounded domain $\Omega$. 
When $\Omega = \mathbb{R}^3$, the lower bound from~\cite{IK1} reads $\|D\omega(t)\|_{L^\infty} \geq C^{-1}(T_\ast -t)^{-7/5}$ for $t$ close to $T_\ast$.
The purpose of the first result of this paper, see Theorem~\ref{T01} and Remark~\ref{R01}, is to show that this lower bound cannot be matched by an upper bound of the same order. 
Consequently, $\|D\omega\|_{L^\infty}$, and more generally $\|D^k\omega\|_{L^\infty}$, cannot blow up exactly at the endpoint rate $(T_\ast -t)^{-1-2k/5}$.

For the second result of this paper, see Theorem~\ref{T03}, we consider the Euler equations with space-time variable coefficients,
given by
  \begin{align}
   \begin{split}
    \partial_t v_i
    + (v_m - \psi_m)a_{km} \partial_k v_i
    + a_{ki} \partial_k q
    &= 0
    \comma 
    i=1,2,3,
    \\
    a_{ji} \partial_j v_i &= 0
    ,
   \end{split}
   \label{EQ02}
  \end{align}
in $\Omega$, with the boundary condition,
  \begin{equation}
   (v_k - \psi_k)a_{jk}n_j = 0
   \inon{on $\partial \Omega$,}
   \label{EQ03}
  \end{equation}
where $\Omega \subset \mathbb{R}^3$ is a smooth bounded domain.
Precise regularity and compatibility properties of the variable coefficients are given in Section~\ref{sec01}. 
The model given by~\eqref{EQ02} and \eqref{EQ03} arises when studying free-boundary problems involving the incompressible Euler equations, 
for which one might fix the boundary via the Arbitrary Lagrangian Eulerian (ALE) change of variables. 
Effectively, the momentum equation~\eqref{EQ02} generalizes the equation obtained from applying a coordinate change on~\eqref{EQ00}. 
In fact, the given function $\psi$ may represent the velocity of the change of variables map which fixes the domain, and its normal component matches the normal velocity of the flow in what is known as the kinematic condition~\eqref{EQ03}. 
For example, the model appears when studying the fluid-structure interaction between an inviscid fluid and an elastic plate, see~\cite{KuT,AKT}.
More generally, the system~\eqref{EQ02}--\eqref{EQ03} models the classical Euler equations in a general time-dependent domain with a kinematic condition which equates the normal fluid velocity with the interface velocity. 
Finally, we note that the model reduces to~\eqref{EQ00} on a bounded domain with the impermeability condition when $a_{ij} = \delta_{ij}$ and $\psi = 0$.

Jointly with Kukavica and Tuffaha~\cite{IKT}, we established the first local-in-time existence result for~\eqref{EQ02} and~\eqref{EQ03} with initial data  $v_0$ in the minimal regularity class $H^r$, where $r > 2.5$.
Furthermore, we derived a BKM-type criterion based on the $\BMO$ norm of the variable vorticity $\zeta_i = \eps_{ijk} b_{\ell j}\partial_\ell v_k$ for solutions in $H^3$. 
In Theorem~\ref{T03}, we expand this blow-up criterion to the minimal regularity class $r \in (2.5,3)$ and give a new criterion involving the $L^\infty$ norm of $\zeta$.

The paper is structured as follows.
In Section~\ref{sec01}, we outline the main results of the paper, Theorems~\ref{T01} and~\ref{T03}, as well as a precise exposition of the system~\eqref{EQ02}--\eqref{EQ03}. 
In Section~\ref{sec02}, we prove Theorem~\ref{T01}, which excludes the endpoint behavior $\|D^k\omega(t)\|_{L^\infty} \sim (T_\ast -t)^{-1-2k/5}$ for every $k \geq 1$. 
The proof employs a result by Chae and Wolf~\cite{CW}, which excludes atomic concentration of the kinetic energy at the blow-up time under a Type~I bound on $\|\nabla u\|_{L^\infty}$.
Assuming that $\|D\omega(t)\|_{L^\infty} \les (T_\ast - t)^{-7/5}$ for late times, a local elliptic estimate yields such a Type~I bound, and the absence of atoms, together with tightness of the kinetic energy, implies that the energy on small balls is uniformly small.
Inserting this back into the local elliptic estimate improves the Type I bound to $\|\nabla u(t)\|_{L^\infty} = o((T_\ast -t)^{-1})$, which contradicts the BKM criterion.
Indeed, the Type I hypothesis needed for using~\cite[Theorem~1.1]{CW} results from the hypothetical endpoint bound itself and, remarkably, requires no additional assumptions on the Euler solution. 
It would be interesting to see if similar results hold for the blowup rate on $\|\omega(t)\|_{L^\infty}$ derived in~\cite{IK1}. 

In Section~\ref{sec03}, we prove Theorem~\ref{T03}, which pertains to the variable Euler system. 
We show that $L_t^1$-integrability of $\|v\|_{H^1} + \|\zeta\|_{\BMO}$ implies continuation of the solution in $H^r$.
The proof employs the Sobolev-Slobodeckij norm, defined in~\eqref{EQ25} and~\eqref{EQ26}, in the evolution estimate of $\|v\|_{H^r}$. 
Finally, we provide a new blow-up criterion involving $\|v\|_{L^2} + \|\zeta\|_{L^\infty}$, which trades the dependence on $\|\nabla v\|_{L^2}$ for the stronger norm $\|\zeta\|_{L^\infty}$.
Regarding the presence of $\|v\|_{L^2}$ in the criterion, a significant trait of the system~\eqref{EQ02}--\eqref{EQ03} is that the $L^2$ norm of the fluid velocity is not conserved. 
In fact, it is not controlled solely by assumed data but is driven by the stretching term, $\|\nabla v\|_{L^\infty}$, due to the inhomogeneous boundary condition.

\section{Preliminaries and main results}\label{sec01}
For results related to Theorem~\ref{T01}, we denote $L^p = L^p(\mathbb{R}^3)$,
and for results regarding Theorem~\ref{T03}, we denote $L^p = L^p(\Omega)$ where $\Omega$ is a smooth bounded domain. 
Also, we use $D$ and $\nabla$ interchangeably to mean the gradient. 

The first theorem pertains to the classical incompressible Euler equations~\eqref{EQ00} posed on $\mathbb{R}^3$.
\begin{Theorem}
\label{T01}
Let $k \geq 1$ be an integer, and let $u$ be a solution to the Euler equations~\eqref{EQ00} on $\mathbb{R}^3$ in the regularity class~\eqref{EQ01}, with $r > k + 5/2$.
Suppose there is a time $T_\ast > 0$ such that the solution cannot be continued in this class to $T = T_\ast$ and that $T_\ast$ is the first such time. Then, we have
  \begin{equation}
    \limsup_{t \to T_\ast^-} (T_\ast - t)^{1+2k/5} \|D^k\omega(t)\|_{L^\infty}
    = 
    \infty
    .
    \label{EQ10}
  \end{equation}
\end{Theorem}
Thus, there is no constant $C$ such that $\|D^k\omega(t)\|_{L^\infty} \leq C(T_\ast -t)^{-1-2k/5}$ for all $t$ sufficiently close to~$T_\ast$. 
\begin{Remark}
\label{R01}
{\rm
Theorem~\ref{T01} complements the pointwise-in-time lower bound established in~\cite{IK1}.
Under the assumptions of Theorem~\ref{T01} with $k=1$, there exists a constant $C$ such that $\|D\omega(t)\|_{L^\infty} \geq C^{-1}(T_\ast -t)^{-7/5}$ for $T_\ast -t < 1/C$, see~\cite{IK1}.
For $k \geq 2$, combining this bound with the interpolation inequality~\eqref{EQ49a} below and the conservation of energy yields $\|D^k\omega(t)\|_{L^\infty} \geq C^{-1}(T_\ast -t)^{-1-2k/5}$ for $T_\ast - t < 1/C$, with a possibly larger constant $C$.
Together with~\eqref{EQ10}, we thus obtain
  \begin{equation}
   0
   <
   \liminf_{t \to T_\ast^-} (T_\ast -t)^{1+2k/5} \|D^k\omega(t)\|_{L^\infty}
   \leq
   \limsup_{t \to T_\ast^-} (T_\ast -t)^{1+2k/5} \|D^k\omega(t)\|_{L^\infty}
   =
   \infty
   \comma
   k \geq 1
   .
   \llabel{EQ10a}
  \end{equation}
Hence, $\|D^k\omega\|_{L^\infty}$ cannot blow up exactly at the endpoint rate $(T_\ast -t)^{-1-2k/5}$: 
it is bounded below by a multiple of this rate for all $t$ close to $T_\ast$, but it is not bounded above by any multiple of it.
Note that Theorem~\ref{T01} and~\cite{IK1} do not determine whether the limit inferior is finite.
}
\end{Remark}

The second result pertains to the variable Euler system given in~\eqref{EQ02}--\eqref{EQ03}.
First, we state some assumptions on the space-time variable coefficients.
Let $\Omega$ be a smooth bounded domain and $r \in (2.5,3)$. 
Let $\psi \colon \Omega \times [0,T] \to \mathbb{R}^3$ and $a \colon \Omega \times [0,T] \to GL_3(\mathbb{R})$ be such that
  \begin{equation}
   (\psi,\psi_t) \in L^\infty([0,T]; H^r \times H^{r-1})
   \label{EQ11}
  \end{equation}
and 
  \begin{equation}
   (a,a_t) \in L^\infty ([0,T]; H^r \times H^{r-1})
   .
   \llabel{EQ12}
  \end{equation}
Next, define $J = \det(a^{-1})$ and $b = \cof(a^{-1})^T$, from which we may write
  \begin{equation}
   a = \frac{1}{J} b
   .
   \llabel{EQ13}
  \end{equation}
We also assume that $J$ is uniformly bounded above and below on $\Omega \times [0,T]$ and that $b$ satisfies the Piola identity,
  \begin{equation}
   \partial_j b_{ji} = 0
   \comma
   i =1,2,3.
   \llabel{EQ14}
  \end{equation}
Finally, we state a compatibility condition which is necessary and sufficient for existence of solutions to~\eqref{EQ02}--\eqref{EQ03}:
  \begin{equation}
   \int_{\partial \Omega} n_j b_{ji} \psi_i = 0
   \comma 
   t \in [0,T]
   .
   \label{EQ15}
  \end{equation}
For details on the necessity of~\eqref{EQ15}, see Remark~4.1 in~\cite{IKT}.
We also note that the pressure $q$ solves the elliptic Neumann problem,
  \begin{equation}
   -\partial_j(b_{ji} a_{ki} \partial_k q)
   =
   -\partial_j(\partial_t b_{ji} v_i)
   + b_{ji} \partial_j((v_m - \psi_m) a_{km})\partial_k v_i
   - (v_m - \psi_m)a_{km} \partial_k b_{ji} \partial_j v_i
   \inon{in $\Omega$}
   \label{EQ15a}
  \end{equation}
and
  \begin{align}
   \begin{split}
    -n_jb_{ji}a_{ki} \partial_k q
    &=
    \partial_t(n_jb_{ji}\psi_i)
    -n_j\partial_t b_{ji} v_i
    +(v_m - \psi_m)a_{km} \partial_k^\tau (n_j b_{ji}\psi_i)
    \\&\indeq\indeq
    - (v_m - \psi_m) a_{km} \partial_k^\tau(n_j b_{ji})v_i
    ,
   \end{split}
   \inon{on $\partial\Omega$,}
   \label{EQ15b}
  \end{align}
where we denote $\partial_k^\tau = (\delta_{k\ell} - n_kn_\ell)\partial_\ell$.
The derivation of~\eqref{EQ15a} and~\eqref{EQ15b} can be found in Section 3.1 of~\cite{IKT}.
Note that the system is of the form $-\Delta q = \partial v \cdot \partial v$ and $\partial_n q = v \cdot v$, which mirrors the system in Lemma~1.1 of~\cite{T}. 

Under the assumptions~\eqref{EQ11}--\eqref{EQ15}, the system given in~\eqref{EQ02}--\eqref{EQ03} admits the following local-in-time existence theorem. 
\begin{Theorem}
\label{T02}
Assume that $v_0 \in H^r$, where $r >2.5$, satisfies the divergence and boundary conditions stated above at time $t =0$.
Let $M > 0$ be such that $\|v_0\|_{H^r} \leq M$. 
Further assume the compatibility condition~\eqref{EQ15}.
Then there exists a local-in-time solution $(v,q)$ to the system~\eqref{EQ02}--\eqref{EQ03} such that
  \begin{align}
   \begin{split}
    v &\in L^\infty([0,T];H^r(\Omega)) \cap C([0,T];H^{r-1}(\Omega)),
    \\
    v_t &\in L^\infty([0,T]; H^{r-1}(\Omega)),
    \\
    q &\in L^\infty([0,T]; H^r(\Omega)),
   \end{split}
   \label{EQ20}
  \end{align}
for some time $T> 0$ depending on $M$ and the assumed data $\psi$ and $a$. 
\end{Theorem}
\noindent
For the proof of Theorem~\ref{T02}, see~\cite[Theorem~2.2]{IKT}. 

The second result of this paper is an improvement on a blow-up criterion for solutions to~\eqref{EQ02}--\eqref{EQ03}
which was originally derived for $r = 3$, see~\cite[Theorem~2.3]{IKT}. 
In particular, we expand the result to the fractional case $r \in (2.5,3)$.
Suppose that
  \begin{equation}
   (\psi,\psi_t) \in L^\infty([0,\infty); H^{r+1} \times H^r)
   \label{EQ21}
  \end{equation}
and 
  \begin{equation}
   (a,a_t) \in L^\infty([0,\infty); H^{r+1}\times H^r)
   .
   \label{EQ22}
  \end{equation}
\begin{Theorem}
\label{T03}
Assume the regularity in~\eqref{EQ21} and~\eqref{EQ22} for $a$ and $\psi$. 
Let $v$ be a solution to~\eqref{EQ02}--\eqref{EQ03} in the class $C([0,T];H^r(\Omega))$ where $r \in (2.5,3)$.
If $T = T_\ast$ is the first time such that $v$ cannot be continued in this class,
then
  \begin{equation}
   \int_0^{T_\ast} \left( \|v(t)\|_{H^1} + \|\zeta(t)\|_{\BMO}\right) \, dt = \infty
   .
   \label{EQ23}
  \end{equation}
Moreover,
  \begin{equation}
          \int_0^{T_\ast} \left( \|v(t)\|_{L^2} + \|\zeta(t)\|_{L^\infty}\right) \, dt = \infty 
   .
   \label{EQ24}
  \end{equation}
\end{Theorem}
\begin{Remark}
{\rm
In fact, the solutions of Theorem~\ref{T02} under the regularity assumptions~\eqref{EQ21}--\eqref{EQ22} 
belong to $C([0,T_\ast);H^r)$.
Indeed, a solution $v$ satisfying~\eqref{EQ20} is in $C_w([0,T];H^r)$.
It then suffices to show that $\|v(t)\|_{H^r}$ is continuous. 
To see this, the estimate~\eqref{EQ51} in conjunction with the approximation argument in Theorem~5.1 in~\cite{IKT} with $v_0$ replaced by $v(s)$ implies
  \begin{equation}
   \|v(t)\|_{H^r}
   \leq 
   \|v(s)\|_{H^r} \exp \left( C\int_s^t \left( 1 +  \|v(\tau)\|_{W^{1,\infty}}\right)\, d\tau\right)
   \label{EQ68}
  \end{equation}
for a solution $v \in L^\infty([0,T];H^r)$. 
Thus, lower and upper semicontinuity of $\|v(t)\|_{H^r}$ follow from weak continuity and~\eqref{EQ68}, respectively.
}
\end{Remark}

\section{Classical Euler: Exclusion of the endpoint rate for $\|D^k\omega\|_{L^\infty}$}\label{sec02}
In this section, all norms are taken over $\mathbb{R}^3$ unless indicated otherwise, and $B_r(x_0)$ denotes the ball of radius $r$ centered at $x_0$.
The next lemma obtains pointwise-in-space bounds on $u$ and $\nabla u$ in terms of the local $L^2$ norm of $u$ and the local maximum norm of $D\omega$.
This lemma will be employed twice in the proof of Theorem~\ref{T01}. 
Under the assumption that $\|D\omega(t)\|_{L^\infty} \les (T_\ast -t)^{-7/5}$, we use it to establish a Type I gradient bound and $L_t^1$ integrability of $\|u(t)\|_{L^\infty}$, which together verify the hypotheses of~\cite[Theorem~1.1]{CW} and Lemma~\ref{L02}. 
Then, critically, we apply it in~\eqref{EQ47} in conjunction with the smallness in~\eqref{EQ46} obtained by Lemma~\ref{L02} to establish the rate $\|\nabla u(t)\|_{L^\infty} = o((T_\ast -t)^{-1})$.
\begin{Lemma}
\label{L01}
Let $u \in C^2(\mathbb{R}^3)$ be a divergence-free vector field and $\omega = \curl u$.
Then, for every $x_0 \in \mathbb{R}^3$ and $r > 0$, we have
  \begin{align}
   |u(x_0)|
   &\les
   r^2 \|D\omega\|_{L^\infty(B_r(x_0))}
   + r^{-3/2} \|u\|_{L^2(B_r(x_0))}
   ,
   \label{EQ30}
   \\
   |\nabla u(x_0)|
   &\les
   r \|D\omega\|_{L^\infty(B_r(x_0))}
   + r^{-5/2} \|u\|_{L^2(B_r(x_0))}
   ,
   \label{EQ31}
  \end{align}
where the implicit constants are absolute.
Consequently, if in addition $u \in L^2(\mathbb{R}^3)$ and $D\omega \in L^\infty(\mathbb{R}^3)$, then
  \begin{equation}
   \|u\|_{L^\infty}
   \les
   \|u\|_{L^2}^{4/7} \|D\omega\|_{L^\infty}^{3/7}
   \andand
   \|\nabla u\|_{L^\infty}
   \les
   \|u\|_{L^2}^{2/7} \|D\omega\|_{L^\infty}^{5/7}
   .
   \label{EQ32}
  \end{equation}
\end{Lemma}
\begin{proof}
By translation, it suffices to consider $x_0 = 0$.
First, for $v \in C^2(\overline{B_1(0)})$, note that interior $W^{2,4}$ estimates for the Laplacian on nested balls give
  \begin{equation}
   \|v\|_{W^{2,4}(B_{1/2}(0))} 
   \les
   \|v\|_{L^4(B_{3/4}(0))}
   + \|\Delta v\|_{L^4(B_{3/4}(0))}
   .
   \llabel{EQ32c}
  \end{equation}
Then, by interior $H^2$ estimates, we have $\|v\|_{H^2(B_{3/4}(0))} \les \|v\|_{L^2(B_1(0))} + \|\Delta v\|_{L^2(B_1(0))}$.
Since $H^2(B) \subset L^4(B)$ for a bounded ball and $W^{2,4} \subset C^1$, we obtain
  \begin{equation}
   |v(0)| + |\nabla v(0)|
   \les
   \|v\|_{W^{2,4}(B_{1/2}(0))}
   \les
   \|v\|_{L^2(B_1(0))}
   + \|\Delta v\|_{L^\infty(B_1(0))}
   .
   \label{EQ32a}
  \end{equation}
Next, since $\div u = 0$, we have $\Delta u = -\curl \omega$.
Applying~\eqref{EQ32a} to $v(y) = u(ry)$, for which $\|v\|_{L^2(B_1(0))} = r^{-3/2}\|u\|_{L^2(B_r(0))}$ and $\|\Delta v\|_{L^\infty(B_1(0))} \les r^2 \|D\omega\|_{L^\infty(B_r(0))}$, we arrive at 
  \begin{equation}
   |u(0)| + r|\nabla u(0)|
   \les
   r^{-3/2} \|u\|_{L^2(B_r(0))}
   + r^2 \|D\omega\|_{L^\infty(B_r(0))}
   ,
   \llabel{EQ32b}
  \end{equation}
which yields~\eqref{EQ30} and~\eqref{EQ31}.

Finally, to get~\eqref{EQ32}, we bound the local $L^2$ norms in~\eqref{EQ30} and~\eqref{EQ31} by $\|u\|_{L^2}$, take the supremum over $x_0 \in \mathbb{R}^3$, and choose $r = (\|u\|_{L^2}/\|D\omega\|_{L^\infty})^{2/7}$.
The case $\|u\|_{L^2} = 0$ is trivial, while if $\|D\omega\|_{L^\infty} = 0$, then $u$ is a square-integrable harmonic function on $\mathbb{R}^3$, and hence $u = 0$.
\end{proof}

The next lemma shows that a uniformly tight family of energy densities converging to a measure without atoms cannot concentrate on small balls, uniformly in time.
Assuming that $\|D\omega(t)\|_{L^\infty} \les (T_\ast -t)^{-7/5}$ for late times, we verify the hypotheses of this lemma in the proof of Theorem~\ref{T01} using a result of~\cite{CW} and a tightness estimate. 
Lemma~\ref{L02} then provides us with the uniform smallness of local energy needed to prove $\|\nabla u(t)\|_{L^\infty} = o((T_\ast -t)^{-1})$.
\begin{Lemma}
\label{L02}
Let $t_0 < T$ and $f \in C([t_0,T); L^2(\mathbb{R}^3))$ with $\sup_{t_0 \leq t < T}\|f(t)\|_{L^2} < \infty$, and denote $\mu_t = |f(t)|^2\, dx$.
Assume that $\mu_t$ converges weakly-$\ast$ to a finite measure $\mu$ as $t \to T^-$, i.e., $\int \varphi\, d\mu_t \to \int \varphi \, d\mu$ for every $\varphi \in C_c(\mathbb{R}^3)$.
Further suppose 
that $\mu(\{x\}) = 0$ for every $x \in \mathbb{R}^3$, and 
  \begin{equation}
   \lim_{R \to \infty} \sup_{t_0 \leq t < T} \mu_t(\mathbb{R}^3 \setminus B_R(0))
   =
   0
   .
   \label{EQ33}
  \end{equation}
Then
  \begin{equation}
   \lim_{r \to 0^+} \sup_{t_0 \leq t < T,\, x \in \mathbb{R}^3} \|f(t)\|_{L^2(B_r(x))}
   =
   0
   .
   \label{EQ34}
  \end{equation}
\end{Lemma}
\begin{proof}
The supremum in~\eqref{EQ34} is nondecreasing in $r$, so if~\eqref{EQ34} fails, there exist $\eps > 0$, $r_n \to 0^+$, $t_n \in [t_0,T)$, and $x_n \in \mathbb{R}^3$ such that $\mu_{t_n}(B_{r_n}(x_n)) \geq \eps$ for all $n$.
By~\eqref{EQ33}, the sequence $(x_n)$ is bounded. 
Indeed, assuming the contrary, we may choose $R > 0$ such that  $B_{r_n}(x_n) \subset \mathbb{R}^3 \setminus B_R(0)$ and~$\sup_{t_0\leq t < T} \mu_t(\mathbb{R}^3 \setminus B_R(0)) < \eps$.
Consequently, $\eps \leq \mu_{t_n}(B_{r_n}(x_n)) \leq \mu_{t_n}(\mathbb{R}^3 \setminus B_R(0)) < \eps$, which is a contradiction.
Next, passing to a subsequence, we may assume that $x_n \to x_\ast$ and $t_n \to t_\ast \in [t_0,T]$.
If $t_\ast < T$, then $f(t_n) \to f(t_\ast)$ in $L^2$ by continuity of $f$.
Since $|B_{r_n}(x_n)| \to 0$, the absolute continuity of the integral then implies
  \begin{equation}
   \mu_{t_n}(B_{r_n}(x_n)) 
   = 
   \int_{B_{r_n}(x_n)} |f(t_n)|^2
   \les
   \|f(t_n) - f(t_\ast)\|_{L^2}^2 
   + \int_{B_{r_n}(x_n)} |f(t_\ast)|^2 
   \to 0
   ,
   \llabel{EQ34b}
  \end{equation}
which is a contradiction.
If $t_\ast = T$, we fix $\delta > 0$ and let $\varphi_\delta \in C_c(\mathbb{R}^3)$ be such that $0 \leq \varphi_\delta \leq 1$, $\varphi_\delta = 1$ on $B_\delta(x_\ast)$, and $\supp \varphi_\delta \subset B_{2\delta}(x_\ast)$.
For $n$ sufficiently large, we have $B_{r_n}(x_n) \subset B_\delta(x_\ast)$, and therefore
  \begin{equation}
   \eps
   \leq
   \int \varphi_\delta \, d\mu_{t_n}
   \to
   \int \varphi_\delta \, d\mu
   \leq
   \mu(B_{2\delta}(x_\ast))
   \comma
   n \to \infty
   .
   \llabel{EQ34a}
  \end{equation}
Letting $\delta \to 0$, we obtain $\mu(\{x_\ast\}) \geq \eps$, contradicting the assumption that $\mu$ has no atoms.
\end{proof}

Now, we prove Theorem~\ref{T01}.
\begin{proof}[Proof of Theorem~\ref{T01}]
We first prove the theorem for $k = 1$.
Towards a contradiction, suppose that the limit superior in~\eqref{EQ10} is finite for $k = 1$.
Then there exist $t_0 \in [0,T_\ast)$ and $K > 0$ such that
  \begin{equation}
   \|D\omega(t)\|_{L^\infty}
   \leq
   \frac{K}{(T_\ast - t)^{7/5}}
   \comma
   t \in [t_0,T_\ast)
   .
   \label{EQ35}
  \end{equation}
Under this assumption, we claim that 
  \begin{equation}
   \lim_{t \to T_\ast^-} (T_\ast - t) \|\nabla u(t)\|_{L^\infty}
   =
   0
   .
   \label{EQ36}
  \end{equation}
If~\eqref{EQ36} holds, then maximality of $T_\ast$ is violated.
Indeed, integrating the vorticity equation
  \begin{equation}
   \omega_t
   + u \cdot \nabla \omega
   =
   \omega \cdot \nabla u
   \inon{in $\mathbb{R}^3$}
   \llabel{EQ37}
  \end{equation}
along particle trajectories, we obtain
  \begin{equation}
   \|\omega(t)\|_{L^\infty}
   \leq
   \|\omega(t_1)\|_{L^\infty}
   + C\int_{t_1}^t \|\nabla u(s)\|_{L^\infty} \|\omega(s)\|_{L^\infty}\, ds
   \comma
   0 \leq t_1 \leq t < T_\ast
   ,
   \llabel{EQ38}
  \end{equation}
where $C$ is an absolute constant.
By~\eqref{EQ36}, we may choose $t_1 \in [t_0,T_\ast)$ such that $C(T_\ast - t)\|\nabla u(t)\|_{L^\infty} \leq 1/2$ for $t \in [t_1,T_\ast)$.
Gronwall's inequality then gives
  \begin{equation}
   \|\omega(t)\|_{L^\infty}
   \leq
   \|\omega(t_1)\|_{L^\infty} \exp\left( \int_{t_1}^t \frac{ds}{2(T_\ast - s)}\right)
   =
   \|\omega(t_1)\|_{L^\infty} \left( \frac{T_\ast - t_1}{T_\ast - t}\right)^{1/2}
   \comma
   t \in [t_1,T_\ast)
   .
   \llabel{EQ39}
  \end{equation}
It follows that $\int_0^{T_\ast} \|\omega(t)\|_{L^\infty}\, dt < \infty$.
By the BKM criterion, the solution $u$ may then be continued in the class~\eqref{EQ01} beyond $T_\ast$, contradicting the definition of $T_\ast$.
Hence, it suffices to establish~\eqref{EQ36}.

Now, we verify the hypotheses of Lemma~\ref{L02}.
We first show that $|u(t)|^2\, dx$ converges weakly-$\ast$ to a finite atomless measure as $t \to T_\ast^-$.
By~\eqref{EQ32}, the bound~\eqref{EQ35}, and the conservation of energy, we have
  \begin{equation}
   \|\nabla u(t)\|_{L^\infty}
   \les
   \|u_0\|_{L^2}^{2/7} K^{5/7}(T_\ast - t)^{-1}
   \andand
   \|u(t)\|_{L^\infty}
   \les
   \|u_0\|_{L^2}^{4/7} K^{3/7}(T_\ast - t)^{-3/5}
   \label{EQ40}
  \end{equation}
for $t \in [t_0,T_\ast)$.
In particular,
  \begin{equation}
   \int_{t_0}^{T_\ast} \|u(t)\|_{L^\infty}\, dt
   <
   \infty
   .
   \label{EQ41}
  \end{equation}
The first bound in~\eqref{EQ40} is a Type~I condition on the velocity gradient.
Moreover, since $r > 7/2$, the solution is bounded in $W^{1,\infty}(\mathbb{R}^3)$ on $[t_0,t]$ for every $t < T_\ast$, and its $L^2$ norm is conserved.
Thus, after a translation and rescaling in time, the result of Chae and Wolf~\cite[Theorem~1.1 and Corollary~1.2]{CW} applies on $\mathbb{R}^3 \times (t_0,T_\ast)$.
That is, the measures $|u(t)|^2\, dx$ converge weakly-$\ast$ to a finite measure $\mu$ as $t \to T_\ast^-$, 
  \begin{equation}
   \int \varphi(x) |u(x,t)|^2\, dx
   \to
   \int \varphi\, d\mu
   \comma
   \varphi \in C_c(\mathbb{R}^3)
   ,
   \llabel{EQ42}
  \end{equation}
and $\mu(\{x\}) = 0$ for every $x \in \mathbb{R}^3$.

Next, we show that the energy is uniformly tight on $[t_0,T_\ast)$.
We write the pressure as $p = R_iR_j(u_iu_j)$, where $R_1, R_2, R_3$ are the Riesz transforms.
By $L^2$ boundedness of the Riesz transforms,
  \begin{equation}
   \|p(t)\|_{L^2}
   \les
   \|u(t)\|_{L^\infty} \|u(t)\|_{L^2}
   =
   \|u_0\|_{L^2} \|u(t)\|_{L^\infty}
   .
   \label{EQ43}
  \end{equation}
Taking the inner product of the momentum equation in~\eqref{EQ00} with $2u$ and using the divergence-free condition, we obtain the local energy identity
  \begin{equation}
   \partial_t |u|^2
   + \div \bigl( (|u|^2 + 2p) u\bigr)
   =
   0
   \inon{in $\mathbb{R}^3$.}
   \label{EQ44}
  \end{equation}
For $R > 0$, let $\chi_R \in C^\infty(\mathbb{R}^3)$ be such that $0 \leq \chi_R \leq 1$, $\chi_R = 0$ on $B_R(0)$, $\chi_R = 1$ outside $B_{2R}(0)$, and $|\nabla \chi_R| \les R^{-1}$.
Testing~\eqref{EQ44} with $\chi_R$ on $[t_0,t]$ and using~\eqref{EQ43}, we integrate by parts\footnote
{
Since $\chi_R$ is not compactly supported, integration by parts in \eqref{EQ45} is justified by first testing~\eqref{EQ44} with $\chi_R\psi_\rho$, where $\psi_\rho \in C_c^\infty(\mathbb{R}^3)$ satisfies $0 \leq \psi_\rho \leq 1$, $\psi_\rho = 1$ on $B_\rho(0)$, $\supp \psi_\rho \subset B_{2\rho}(0)$, and $|\nabla \psi_\rho| \les \rho^{-1}$, and then letting $\rho \to \infty$; 
the resulting error term is bounded by a constant multiple of $\rho^{-1}\int_{t_0}^{T_\ast} \int (|u|^3 + |p||u|)\, dx ds$, which is finite by~\eqref{EQ41} and~\eqref{EQ43}.
}
to get
  \begin{equation}
    \int_{|x| > 2R} |u(x,t)|^2\, dx
    \leq
    \int_{|x| > R} |u(x,t_0)|^2 \, dx
    + \frac{C\|u_0\|_{L^2}^2}{R} \int_{t_0}^{T_\ast} \|u(s)\|_{L^\infty}\, ds
    \comma
    t \in [t_0,T_\ast)
    .
   \label{EQ45}
  \end{equation}
By~\eqref{EQ41}, the right-hand side of~\eqref{EQ45} tends to $0$ as $R \to \infty$, uniformly in $t \in [t_0,T_\ast)$.

We may now apply Lemma~\ref{L02} with $f = u$ and $T = T_\ast$, which yields
  \begin{equation}
   m(r)
   :=
   \sup_{t_0 \leq t < T_\ast,\, x \in \mathbb{R}^3} \|u(t)\|_{L^2(B_r(x))}
   \to 0
   \comma
   r \to 0^+
   .
   \label{EQ46}
  \end{equation}
Bounding the local $L^2$ norm in~\eqref{EQ31} by $m(r)$ and taking the supremum over $x_0 \in \mathbb{R}^3$, we obtain
  \begin{equation}
   \|\nabla u(t)\|_{L^\infty}
   \les
   r \|D\omega(t)\|_{L^\infty}
   + r^{-5/2} m(r)
   \comma
   r > 0,\ t \in [t_0,T_\ast)
   .
   \label{EQ47}
  \end{equation}
Now, fix $\lambda > 0$ and set $r = \lambda (T_\ast - t)^{2/5}$ in~\eqref{EQ47}.
By~\eqref{EQ35}, we get
  \begin{equation}
   (T_\ast - t) \|\nabla u(t)\|_{L^\infty}
   \les
   \lambda K
   + \lambda^{-5/2} m\bigl( \lambda (T_\ast - t)^{2/5}\bigr)
   \comma
   t \in [t_0,T_\ast)
   ,
   \llabel{EQ48}
  \end{equation}
where the implicit constant does not depend on $\lambda$.
Letting $t \to T_\ast^-$ and using~\eqref{EQ46}, we arrive at
  \begin{equation}
   \limsup_{t \to T_\ast^-} (T_\ast - t) \|\nabla u(t)\|_{L^\infty}
   \les
   \lambda K
   .
   \llabel{EQ49}
  \end{equation}
Since $\lambda > 0$ is arbitrary, we obtain~\eqref{EQ36}, which completes the proof for $k = 1$.

It remains to consider $k \geq 2$.
We use the interpolation inequality
  \begin{equation}
   \|D\omega\|_{L^\infty}
   \les
   \|u\|_{L^2}^{2(k-1)/(2k+5)} \|D^k\omega\|_{L^\infty}^{7/(2k+5)}
   ,
   \label{EQ49a}
  \end{equation}
where the implicit constant depends on $k$.
Indeed, denoting by $\Delta_j$ the homogeneous Littlewood-Paley projections, Bernstein's inequality gives
  \begin{equation}
   \sum_{j \leq J} \|\Delta_j D\omega\|_{L^\infty}
   \les
   2^{7J/2} \|u\|_{L^2}
   \andand
   \sum_{j > J} \|\Delta_j D\omega\|_{L^\infty}
   \les
   2^{-(k-1)J} \|D^k\omega\|_{L^\infty}
   \comma
   J \in \mathbb{Z}
   ,
   \llabel{EQ49b}
  \end{equation}
and~\eqref{EQ49a} follows by choosing $2^J$ comparable to $(\|D^k\omega\|_{L^\infty}/\|u\|_{L^2})^{2/(2k+5)}$.
Now, suppose that the limit superior in~\eqref{EQ10} is finite for some $k \geq 2$.
Then there exist $t_0 \in [0,T_\ast)$ and $K > 0$ such that $\|D^k\omega(t)\|_{L^\infty} \leq K(T_\ast -t)^{-(2k+5)/5}$ for $t \in [t_0,T_\ast)$.
Since $\|u(t)\|_{L^2} = \|u_0\|_{L^2}$, the inequality~\eqref{EQ49a} gives
  \begin{equation}
   \|D\omega(t)\|_{L^\infty}
   \les
   \frac{\|u_0\|_{L^2}^{2(k-1)/(2k+5)} K^{7/(2k+5)}}{(T_\ast - t)^{7/5}}
   \comma
   t \in [t_0,T_\ast)
   .
   \llabel{EQ49c}
  \end{equation}
As $r > k + 5/2 > 7/2$, this contradicts the case $k = 1$, and the proof is complete.
\end{proof}

\section{Variable Euler: A BKM-type criterion for $r \in (2.5,3)$}\label{sec03}
On a smooth bounded domain, we may identify the fractional Sobolev spaces with Sobolev-Slobodeckij spaces.
For $p \in (1,\infty)$ and $s \in (0,1)$, the Slobodeckij seminorm is defined as 
  \begin{equation}
   [f]_{s,p}
   =
   \left( \int_\Omega \int_\Omega \frac{|f(x) - f(y)|^p}{|x-y|^{3 + sp}} \, dx dy\right)^{1/p}
   ,
   \label{EQ25}
  \end{equation}
  and, for $m \in \{0,1,\dots\}$, the full $W^{m+s,p}$ norm is defined as 
  \begin{equation}
   \|f\|_{W^{m+s,p}} 
   =
   \|f\|_{W^{m, p}}
   + \sup_{|\alpha| = m} [D^\alpha f]_{s,p}
   .
   \label{EQ26}
  \end{equation}
For the rest of the paper, all integrals are over $\Omega$ unless indicated otherwise.

In this section, we prove Theorem~\ref{T03}. 
Throughout the proof, we denote ``$\les$'' to mean ``$\leq$'' with an implicit constant which may depend on $\|\psi\|_{H^{r+1}}$ and $\|a\|_{H^{r+1}}$.
\begin{proof}[Proof of Theorem~\ref{T03}]
Towards a contradiction, suppose that 
  \begin{equation}
   \int_0^{T_\ast} (\|v(t)\|_{H^1} + \|\zeta(t)\|_{\BMO}) \, dt < \infty
   .
   \label{EQ50}
  \end{equation}
We claim that 
  \begin{equation}
   \frac{d}{dt} \|v\|_{H^r}
   \les 
   (1+\|v\|_{W^{1,\infty}})\|v\|_{H^r}
   .
   \label{EQ51}
  \end{equation}
In this proof, we derive~\eqref{EQ51} for a smoother solution $v \in C([0,T];H^{r+1}(\Omega) )$; 
the general case follows by approximation.
The details of this approximation and the subsequent derivation of a contradiction to~\eqref{EQ50} are identical to those in~\cite[Theorem~2.3]{IKT}, and thus we omit the details here.
Hence, it remains only to establish~\eqref{EQ51}. 

We proceed with a symbolic proof of~\eqref{EQ51}.
Applying $D^2$ to~\eqref{EQ02}, we get 
  \begin{equation}
   \partial_t D^2 v_i
   + (v_m - \psi_m)a_{km} \partial_k D^2 v_i
   + D^2(a_{ki} \partial_k q)
   = S_i
   \comma
   i=1,2,3,
   \label{EQ52}
  \end{equation}
where we define
  \begin{equation}
   S_i
   = 
   (v_m - \psi_m)a_{km} \partial_k D^2 v_i -D^2((v_m-\psi_m) a_{km} \partial_k v_i)
   .
   \llabel{EQ53}
  \end{equation}
We consider~\eqref{EQ52} evaluated at $x,y \in \Omega$, with  $x \neq y$, and subtract the two equations.
Denoting $\delta f = f(x) - f(y)$ and $\phi_k = (v_m - \psi_m)a_{km}$, this gives 
  \begin{equation}
   \partial_t \delta (D^2v_i)
   + \phi_k(x) \partial_{x_k}\delta (D^2 v_i)
   + \phi_k(y) \partial_{y_k}\delta (D^2 v_i)
   + Q_i(x) - Q_i(y)
   =
   S_i(x) - S_i(y)
   ,
   \llabel{EQ54}
  \end{equation}
where we define $Q_i(x) = D_x^2(a_{ki}(x)\partial_{x_k}q(x)) $. 
Multiplying this equation by $\delta (D^2v_i) K(x,y)$, where $K(x,y) = |x-y|^{-3-2s}$ is the kernel for the $\dot{W}^{s,2}$ seminorm and $s := r-2 \in (.5,1)$, leads to
  \begin{align}
   \begin{split}
    &\frac{1}{2}\frac{d}{dt} \iint |\delta (D^2 v)|^2 K(x,y) \, dx dy
    \\&\indeq\indeq=
    -\iint \phi_k(x) \partial_{x_k} \delta (D^2 v_i) \delta (D^2 v_i) K(x,y)\, dx dy
    - \iint \phi_k(y) \partial_{y_k} \delta(D^2 v_i) \delta (D^2 v_i) K(x,y) \, dxdy
    \\&\indeq\indeq\indeq
    - \iint (Q_i(x) - Q_i(y)) \delta(D^2 v_i) K(x,y)\, dx dy
    + \iint (S_i(x) - S_i(y)) \delta(D^2 v_i) K(x,y)\, dx dy
    \\&\indeq\indeq
    =: I_1 + I_2 + I_3 + I_4
    .
   \end{split}
   \label{EQ55}
  \end{align}
Integration by parts\footnote{Here, $I_1 + I_2$ is understood as a single integral, whose integrand $-\delta(\phi_k \partial_k D^2 v_i)\, \delta(D^2 v_i) K(x,y)$ is integrable for $v \in H^{3+s}$. To justify~\eqref{EQ55a}, we replace $K$ by the bounded kernel $K_\eps(x,y) = (|x-y|^{3+2s} + \eps)^{-1}$, integrate by parts, where the boundary terms vanish since $\phi_k n_k = 0$ on $\partial\Omega$ by~\eqref{EQ03}, and let $\eps \to 0$ using the dominated convergence theorem, noting that $K_\eps \leq K$ and $|\nabla_x K_\eps| \leq |\nabla_x K|$.} on the first two terms gives 
  \begin{align}
   \begin{split}
    I_1 + I_2
    &=
    \frac{1}{2}\iint (\partial_{x_k} \phi_k(x) + \partial_{y_k} \phi_k(y)) |\delta(D^2 v)|^2 K(x,y)\, dx dy
    \\&\indeq
    + \frac{1}{2}\iint (\phi_k(x)- \phi_k(y)) |\delta (D^2 v)|^2 \partial_{x_k} K(x,y)\, dx dy
    .
   \end{split}
   \label{EQ55a}
  \end{align}
Therefore, we have
  \begin{equation}
   |I_1 + I_2| 
   \les
   \|D\phi\|_{L^\infty} [D^2 v]_{s,2}^2
   \les 
   (1 + \|v\|_{W^{1,\infty}}) [D^2 v]_{s,2}^2
   . 
   \label{EQ56}
  \end{equation}
Next, H\"older's inequality yields
  \begin{equation}
   |I_3|
   \leq
   [Q]_{s,2} [D^2v]_{s,2}
   \les
   \|\nabla q\|_{H^{r}}[D^2v]_{s,2}
   \les
   (1 + \|v\|_{W^{1,\infty}})\|v\|_{H^r}^2
   ,
   \llabel{EQ57}
  \end{equation}
where, in the last inequality, we used elliptic regularity for $q$.
Indeed, under the higher regularity assumptions~\eqref{EQ21} and~\eqref{EQ22}, the elliptic Neumann problem~\eqref{EQ15a}--\eqref{EQ15b} implies 
  \begin{equation}
   \|\nabla q\|_{H^r}
   \les
   (1+\|v\|_{W^{1,\infty}})\|v\|_{H^r}
   .
   \llabel{EQ57a}
  \end{equation}
For $I_4$, we have 
  \begin{equation}
   |I_4|
   \leq
   [S]_{s,2}[D^2v]_{s,2}
   .
   \label{EQ58}
  \end{equation}
To estimate $[S]_{s,2}$, we denote by $\,\tilde \, :H^{2+s}(\Omega) \to H^{2+s}(\mathbb{R}^3)$ a continuous extension operator which is also continuous from $W^{1,\infty}(\Omega) \to W^{1,\infty}(\mathbb{R}^3)$. 
We estimate 
  \begin{align}
   \begin{split}
    [S_i]_{s,2}
    &=
    [D^2(\phi_k\partial_k v_i) - \phi_k \partial_k D^2 v_i]_{s,2}
    \leq
    [D^2(\tilde{\phi}_k \partial_k \tilde{v}_i) - \tilde{\phi}_k \partial_k D^2 \tilde{v}_i]_{s,2,\mathbb{R}^3}
    \\&\les
    \|D^2 (\tilde{\phi}_k \partial_k \tilde{v}_i) - \tilde{\phi}_k \partial_k D^2 \tilde{v}_i\|_{H^s(\mathbb{R}^3)}
    \leq
    \|D^2 \tilde{\phi}_k D\tilde{v}_i\|_{H^s(\mathbb{R}^3)}
    + \|D\tilde{\phi}_k D^2 \tilde{v}_i\|_{H^s(\mathbb{R}^3)}
    .
   \llabel{EQ59}
   \end{split}
  \end{align}
We consider the first term on the right-hand side as the other term is bounded similarly. 
By the fractional Leibniz rule, 
  \begin{equation}
   \|D^2 \tilde{\phi}_k D\tilde{v}_i\|_{H^s(\mathbb{R}^3)}
   \les
   \|D^{2}\tilde{\phi}\|_{H^s(\mathbb{R}^3)} \|D\tilde{v}\|_{L^\infty(\mathbb{R}^3)}
   + \|D^2 \tilde{\phi}\|_{L^p(\mathbb{R}^3)}\|D\tilde{v}\|_{W^{s,q}(\mathbb{R}^3)}
   ,
   \llabel{EQ60}
  \end{equation}
where $p,q \in (1,\infty)$ satisfy $p^{-1} + q^{-1} = 1/2$.
We fix $p \in (3, 6(3-2s)^{-1})$, and hence $q \in (3s^{-1}, 6)$.
Then, by the Gagliardo-Nirenberg inequalities,
  \begin{align}
   \begin{split}
    \|D^2 \tilde{\phi}\|_{L^p(\mathbb{R}^3)}
    &\les
    \|\tilde{\phi}\|_{W^{1,\infty}(\mathbb{R}^3)}^{1-\theta_1}\|\tilde{\phi}\|_{H^{2+s}(\mathbb{R}^3)}^{\theta_1}
    \comma 
    \text{where }
    \theta_1 = \frac{1-3/p}{s - 1/2}
    ,
    \\
    \|D\tilde{v}\|_{W^{s,q}(\mathbb{R}^3)}
    &\les
    \|\tilde{v}\|_{W^{1,\infty}(\mathbb{R}^3)}^{1-\theta_2}\|\tilde{v}\|_{H^{2+s}(\mathbb{R}^3)}^{\theta_2}
    \comma
    \text{where }
    \theta_2 = \frac{s - 3/q}{s- 1/2}
    .
   \end{split}
   \llabel{EQ61}
  \end{align}  
Since $\theta_1 + \theta_2 = 1$, we have
  \begin{equation}
   [S_i]_{s,2}
   \les
   (1+\|v\|_{W^{1,\infty}}) \|v\|_{H^r}
   \label{EQ62}
  \end{equation}
by continuity of the extension.
Upon combining~\eqref{EQ56}--\eqref{EQ58} and~\eqref{EQ62}, we arrive at 
  \begin{equation}
   \frac{d}{dt} [D^2v]_{s,2}
   \les 
   (1+\|v\|_{W^{1,\infty}}) \|v\|_{H^r}
   .
   \label{EQ63}
  \end{equation}
As in~\cite[Theorem~2.3]{IKT}, we also have
  \begin{equation}
   \frac{d}{dt} \|v\|_{H^2}
   \les
   \left( 1 + \|v\|_{W^{1,\infty}} \right) \|v\|_{H^2}
   .
   \llabel{EQ64}
  \end{equation}
Combining this with~\eqref{EQ63}, we obtain~\eqref{EQ51}. 

Now, we show that~\eqref{EQ24} holds.
By Lemma~3.5 in~\cite{IKT}, we know that
  \begin{equation}
   \|v\|_{H^1}
   \les
   \|\zeta\|_{L^2}
   +\|n_jb_{jk}\psi_k\|_{H^{1/2}(\partial \Omega)}
   + \|v\|_{L^2}
   \les
   \|\zeta\|_{L^\infty}
   + \|a\|_{H^r}\|\psi\|_{H^r}
   + \|v\|_{L^2}
   ,
   \label{EQ66}
  \end{equation}
where in the second inequality we used that $\Omega$ is bounded. 
Since $a,\psi \in L^\infty([0,\infty) ; H^{r+1}(\Omega))$, the estimate~\eqref{EQ66}, in conjunction with~\eqref{EQ23}, implies 
  \begin{equation}
   \int_0^{T_\ast} \left(\|v(t)\|_{L^2} + \|\zeta(t)\|_{L^\infty}\right)\,dt
   =
   \infty
   ,
   \llabel{EQ67}
  \end{equation}
where we bounded $\|\zeta\|_{\BMO} \leq \|\zeta\|_{L^\infty}$.
This completes the proof of~\eqref{EQ24}. 
\end{proof}

\section*{Acknowledgments}
The author was supported in part by the NSF grant DMS-2205493.      
The author used ChatGPT 6.1 Sol for writing and some arguments in Section~\ref{sec02}. 
No AI was used in the preparation of Section~\ref{sec03}.


\begin{thebibliography}{[BGLT2]}
\small
\bibitem[A]{A}
  M.S.~Aydin,
  \textit{A BKM-type criterion for the Euler equations},
  arXiv:2505.18304.

\bibitem[AKT]{AKT}
M.~Aydin, I.~Kukavica, and A.~Tuffaha,
  \emph{A free boundary inviscid model of flow-structure interaction},
  arXiv:2501.12515.

\bibitem[BKM]{BKM} %MR763762
 J.T.~Beale, T.~Kato, and A.~Majda,
  \emph{Remarks on the breakdown of smooth solutions for the {$3$}-{D} {E}uler equations},
  Comm. Math. Phys.~\textbf{94} (1984), no.~1, 61--66.
  %\MR{763762}

\bibitem[C]{C} %MR2134082
D.~Chae,
  \emph{Remarks on the blow-up criterion of the three-dimensional {E}uler equations},
  Nonlinearity~\textbf{18} (2005), no.~3, 1021--1029.
  %\MR{2134082}

\bibitem[Ch]{Ch}
J.~Chen,
  \emph{Asymptotically self-similar blowup for 3D incompressible Euler with $C^{1,1/3-}$ velocity II: 3D profiles, blowup, and limiting behavior},
  arXiv:2605.15130.

\bibitem[CFM]{CFM} %MR1387460
  P.~Constantin, C.~Fefferman, and A.J.~Majda,
  \emph{Geometric constraints on potentially singular solutions for the {$3$}-{D} {E}uler equations},
  Comm. Partial Differential Equations~\textbf{21} (1996), no.~3-4, 559--571.
    %\MR{1387460}

\bibitem[CIV]{CIV}
P.~Constantin, M.~Ignatova, and V.~Vicol,
  \emph{On putative self-similarity for incompressible 3D Euler},
  arXiv:2602.17570.

\bibitem[CMZ]{CMZ} %MR4929621
D.~C\'ordoba, L.~Martinez-Zoroa, and F.~Zheng,
  \emph{Finite time singularities to the 3{D} incompressible {E}uler
  equations for solutions in {$C^\infty (\Bbb R^3 \setminus \{0\})\cap C^{1,\alpha}\cap L^2$}},
  Ann. PDE~\textbf{11} (2025), no.~2, Paper No. 19, 56.
    %\MR{4929621}

\bibitem[CP]{CP} %MR2954517
T.~Chen and N.~Pavlovi\'c,
  \emph{A lower bound on blowup rates for
  the 3{D} incompressible {E}uler equation and a single exponential
  {B}eale-{K}ato-{M}ajda type estimate}, Comm. Math. Phys.~\textbf{314} (2012),
  no.~1, 265--280.
    %\MR{2954517}

\bibitem[CW]{CW}
D.~Chae and J.~Wolf,
  \emph{Energy concentrations and {T}ype {I} blow-up for the 3{D} {E}uler equations},
  arXiv:1706.02020.

\bibitem[E]{E} %MR4334974
T.~Elgindi,
  \emph{Finite-time singularity formation for {$C^{1,\alpha}$} solutions to the incompressible {E}uler equations on {$\Bbb R^3$}},
  Ann.~of Math.~(2)~\textbf{194} (2021), no.~3, 647--727.
   %\MR{4334974}

\bibitem[F]{F} %MR1230029
A.B.~Ferrari,
  \emph{On the blow-up of solutions of the {$3$}-{D} {E}uler equations in a bounded domain},
  Comm. Math. Phys.~\textbf{155} (1993), no.~2, 277--294.
  %\MR{1230029}

\bibitem[IK1]{IK1}
B.~Ingimarson, I.~Kukavica,
  \textit{Lower bounds on the blowup rate of vorticity in the Euler equations},
  arXiv:2603.17431.

\bibitem[IK2]{IK2}
B.~Ingimarson, I.~Kukavica,
  \emph{On the blowup rate of vorticity for the Euler equations in a bounded domain},
  arXiv:2604.21299.

\bibitem[IKT]{IKT}
B.~Ingimarson, I.~Kukavica, and A.~Tuffaha,
  \emph{The Euler equations with variable coefficients},
  arXiv:2509.01067.

\bibitem[JMO]{JMO}
I.-J.~Jeong, L.~Martínez-Zoroa, and W.S.~O\.za\'nski,
  \emph{Instantaneous continuous loss of Sobolev regularity for the 3D incompressible Euler equation},
  arXiv:2508.06333.

\bibitem[KOT]{KOT} %MR1980623
H.~Kozono, T.~Ogawa, and Y.~Taniuchi,
  \emph{The critical {S}obolev inequalities in {B}esov spaces and regularity criterion to some
  semi-linear evolution equations},
  Math. Z.~\textbf{242} (2002), no.~2, 251--278.
    %\MR{1980623}

\bibitem[KT]{KT} %MR1794270
H.~Kozono and Y.~Taniuchi,
  \emph{Limiting case of the {S}obolev inequality in {BMO}, with application to the {E}uler equations},
  Comm. Math. Phys.~\textbf{214} (2000), no.~1, 191--200.
    %\MR{1794270}

\bibitem[KuT]{KuT}
I.~Kukavica and A.~Tuffaha,
  \emph{A free boundary inviscid model of flow-structure interaction},
  J. Differential Equations 413 (2024), 851-912.

\bibitem[MB]{MB} %MR1867882
A.J.~Majda and A.L.~Bertozzi,
  \emph{Vorticity and incompressible flow},
  Cambridge Texts in Applied Mathematics, vol.~27, Cambridge University Press, Cambridge, 2002.
  %\MR{1867882}


\bibitem[MP]{MP} %MR1245492
C.~Marchioro and M.~Pulvirenti,
  \emph{Mathematical theory of incompressible nonviscous fluids},
  Applied Mathematical Sciences, vol.~96, Springer-Verlag, New York,~1994.
    %\MR{1245492}

\bibitem[P]{P} %MR788778
G.~Ponce,
  \emph{Remarks on a paper: ``{R}emarks on the breakdown of smooth solutions for the {$3$}-{D} {E}uler equations'' by {J}.~{T}.~{B}eale, {T}.~{K}ato and {A}.~{M}ajda}, Comm. Math. Phys.~\textbf{98} (1985), no.~3,
  349--353.
    %\MR{788778}

\bibitem[S]{S}
S.~Shkoller,
  \emph{Incompressible Euler blowup at the $C^{1,\frac13}$ threshold},
  arXiv:2603.10945.

\bibitem[SWZZ]{SWZZ}
F.~Shao, D.~Wei, P.~Zhang, and Z.~Zhang,
  \emph{Self-similar blow-up solutions of incompressible Euler equations in $\mathbb{R}^d$, $d \geq 3$ with $C^{1,1-2/d-}$ velocity},
  arXiv:2605.19716.

\bibitem[SY]{SY} %MR1222827
T.~Shirota and T.~Yanagisawa,
  \emph{A continuation principle for the {$3$}-{D} {E}uler equations for incompressible fluids in a bounded domain},
  Proc.\ Japan Acad.\ Ser.~A Math.\ Sci.~\textbf{69} (1993), no.~3, 77--82.
  %\MR{1222827}

\bibitem[T]{T} %MR393898
R.~Temam,
  \emph{On the {E}uler equations of incompressible perfect fluids},
  S\'eminaire {G}oulaouic-{L}ions-{S}chwartz 1974--1975: \'Equations aux
  d\'eriv\'ees partielles lin\'eaires et non lin\'eaires, \'Ecole Polytech.,
  Paris, 1975, pp.~Exp. No. 10, 15.
    %\MR{393898}

\bibitem[Z]{Z} %MR1101466
W.M.~Zajaczkowski,
  \emph{Remarks on the breakdown of smooth solutions for the {$3$}-d {E}uler equations in a bounded domain},
  Bull.\ Polish Acad.\ Sci.\ Math.~\textbf{37} (1989), no.~1-6, 169--181.
    %\MR{1101466}

\end{thebibliography}
\end{document}